\documentclass[12pt]{amsart}

\usepackage[lite]{amsrefs}

\usepackage{verbatim}
\usepackage{amssymb}
\usepackage{upref}
\usepackage[all]{xy}

\newtheorem{thm}{Theorem}[section]
\newtheorem{cor}[thm]{Corollary}
\newtheorem{lem}[thm]{Lemma}
\newtheorem{prop}[thm]{Proposition}
\newtheorem{obs}[thm]{Observation}

\theoremstyle{definition}
\newtheorem{defn}[thm]{Definition}
\newtheorem{ex}[thm]{Example}

\theoremstyle{remark}
\newtheorem*{rem}{Remark}

\newcommand{\secref}[1]{Section~\textup{\ref{#1}}}
\newcommand{\thmref}[1]{Theorem~\textup{\ref{#1}}}
\newcommand{\corref}[1]{Corollary~\textup{\ref{#1}}}
\newcommand{\lemref}[1]{Lemma~\textup{\ref{#1}}}
\newcommand{\propref}[1]{Proposition~\textup{\ref{#1}}}

\renewcommand{\)}{\textup)}

\newcommand{\righttext}[1]{\qquad\text{#1 }}
\renewcommand{\for}{\righttext{for}}
\renewcommand{\forall}{\righttext{for all}}
\newcommand{\midtext}[1]{\quad\text{ #1 }\quad}
\renewcommand{\and}{\midtext{and}}

\DeclareMathOperator{\aut}{Aut}

\DeclareMathOperator*{\spn}{span}
\DeclareMathOperator*{\clspn}{\overline{\spn}}

\DeclareMathOperator{\ran}{\text{range}}

\newcommand{\N}{\mathbb N}
\newcommand{\Z}{\mathbb Z}

\newcommand{\C}{\mathbb C}

\renewcommand{\H}{\mathcal H}
\newcommand{\G}{\mathcal G}
\newcommand{\K}{\mathcal K}
\renewcommand{\C}{\mathcal C}

\newcommand{\A}{\mathcal A}

\renewcommand{\P}{\mathcal P}

\newcommand{\LL}{\Lambda}
\newcommand{\OO}{\Omega}
\renewcommand{\S}{\Sigma}
\newcommand{\DD}{\Delta}

\renewcommand{\a}{\alpha}
\renewcommand{\b}{\beta}
\newcommand{\e}{\eta}
\newcommand{\ee}{\varepsilon}
\newcommand{\g}{\gamma}
\renewcommand{\d}{\delta}
\renewcommand{\ll}{\lambda}
\newcommand{\m}{\mu}
\newcommand{\n}{\nu}
\newcommand{\p}{\phi}
\renewcommand{\t}{\theta}

\renewcommand{\k}{\kappa}

\newcommand{\Chi}{\raisebox{2pt}{\ensuremath{\chi}}}

\newcommand{\id}{\textnormal{id}}

\renewcommand{\:}{\colon}

\renewcommand{\iff}{\Leftrightarrow}

\newcommand{\inv}{^{-1}}

\newcommand{\what}{\widehat}

\newcommand{\pause}{\renewcommand{\qed}{}\end{proof}}
\newcommand{\resume}[1]{\begin{proof}[Back to the proof of #1]}
\newcommand{\proofof}[1]{\begin{proof}[Proof of #1]}

\newcommand{\fillsquare}[5]{
\xymatrix
@C-5pt
@R-5pt
{
.
\ar@{}[dr]|#5
\ar[d]_#1
&.
\ar[l]_#2
\ar[d]^#3
\\
.
&.
\ar[l]^#4
}
}

\newcommand{\point}{*{\bullet}}

\newcommand{\ff}[1]{\mathcal P(#1)}
\renewcommand{\gg}[1]{\mathcal G(#1)}

\begin{document}

\title[$k$-graph coverings]{Coverings of $k$-graphs}

\author[Pask]{David Pask}
\address{School of Mathematical and Physical Sciences
\\University of Newcastle
\\
NSW 2308
\\
Australia}
\email{davidp@maths.newcastle.edu.au}

\author[Quigg]{John Quigg}
\address{Department of Mathematics and Statistics
\\Arizona State University
\\Tempe, Arizona 85287
\\USA}
\email{quigg@math.asu.edu}

\author[Raeburn]{Iain Raeburn}
\address{School of Mathematical and Physical Sciences
\\University of Newcastle
\\
NSW 2308
\\
Australia}
\email{iain@maths.newcastle.edu.au}

\thanks{This research was supported by grants
from the Australian Research Council and the University of Newcastle}

%



\begin{abstract}
$k$-graphs are higher-rank analogues of directed graphs which were first
developed to provide combinatorial models for operator algebras of
Cuntz-Krieger type. Here we develop the theory of covering spaces for
$k$-graphs, obtaining a satisfactory version of the usual topological
classification in terms of subgroups of a fundamental group. We then use this
classification to describe the $C^*$-algebras of covering $k$-graphs as
crossed products by coactions of homogeneous spaces, generalizing recent
results on the $C^*$-algebras of graphs.
\end{abstract}

\maketitle


\section{Introduction}
\label{sec:intro}

\emph{$k$-graphs} are combinatorial structures which are
$k$-dimensional analogues of (directed) graphs. They were introduced by
Kumjian and the first author
\cite{kp:kgraph} to help understand work of Robertson and
Steger on higher-rank analogues of the Cuntz-Krieger algebras
\cites{rob-steg, rob-steg2}.
The theory of $k$-graphs and their $C^*$-algebras  parallels in many
respects that of graphs and their Cuntz-Krieger algebras
\cites{kp:kgraph, rsy, rsy_finitely}.
Here we investigate to what extent there is an analogue
for $k$-graphs of the theory of coverings of graphs, and the
implications of this theory for the
$C^*$-algebras of $k$-graphs.

A \emph{covering}  of a graph $F$ is by definition a surjective
graph morphism $p\:E\to F$ which is a local isomorphism. As for coverings
of topological spaces, the coverings of $F$ are classified by the
conjugacy classes of subgroups of the fundamental group $\pi_1(F)$, and
every connected covering arises as a quotient of a universal covering
(see \cite{gt} or \cite{dpr}, for example). This last theorem has
interesting ramifications for the Cuntz-Krieger algebras $C^*(E)$ of
covering graphs: if $p\:E\to F$ is a covering, then there is a coaction
$\delta$ of
the fundamental group $\pi_1(F)$ on $C^*(F)$ and a subgroup $H$ of
$\pi_1(F)$ such that $C^*(E)$ is isomorphic to the crossed product
$C^*(F)\times_{\delta} (\pi_1(F)/H)$ \cite{dpr}*{Theorem~3.2}. This
theorem has in turn been of considerable interest in nonabelian
duality for $C^*$-algebras: it provided a family  of crossed
products by homogeneous spaces which we could analyse using our
understanding of graph algebras, and this analysis inspired substantial
improvements in Mansfield's Imprimitivity Theorem \cite{hr:mansfield}.

We seek, therefore,
an analogue of this theory of covering graphs for
$k$-graphs, and a generalisation of \cite{dpr}*{Theorem~3.2} which
describes the $C^*$-algebras of the covering $k$-graphs. Any theory of
coverings must involve the fundamental group, and a majority of the
authors prefer to use the whole fundamental groupoid. We showed in
\cite{pqr:groupoid} that the fundamental groupoids of
$k$-graphs do not behave as well as one might hope, and in particular
that the path category need not embed faithfully in the fundamental
groupoid when $k>1$. So it is something of a relief that our final
results on coverings mirror in every respect the classical topological
theory.

Our approach is to exploit an equivalence  between the coverings of a
$k$-graph and actions of its fundamental groupoid, under which the
connected coverings correspond to transitive actions. Thus we deduce
many of our main theorems from a classification of the transitive
actions of an arbitrary groupoid.

Because every small category is isomorphic to a quotient of a path category,
it will be clear from the proofs that all our results carry over to arbitrary
small categories; however,
we eschew such a generalization since we have no useful applications.

\medskip

After we completed this paper, we learned of the existence of
\cites{bridson, rbrown, higgins, kum:fundgroupoid},
which contain results similar to some of ours.
In \cite{bridson}*{Appendix}, Bridson and Haefliger
develop the elementary theory of
the fundamental group
and coverings
of a small category and proves results similar to some
of ours.
Bridson and Haefliger concentrate on the fundamental \emph{group} ---
indeed, they stop just short of defining the fundamental groupoid.
In
\cites{rbrown, higgins}, Brown and Higgins
investigate coverings of groupoids, and prove the
equivalence with groupoid actions.
Our work was done completely independently of these other sources, and we
believe our methods are of interest, especially our use of skew products.
In
\cite{kum:fundgroupoid}, Kumjian
develops, in the specific
context of $k$-graphs, the fundamental groupoid and the existence of the
universal covering, and proves that, under reasonable hypotheses,
the $C^*$-algebra of the universal covering $k$-graph is Rieffel-Morita
equivalent to a commutative algebra.
We thank Kumjian for bringing \cite{bridson} to our attention.

\medskip

We begin in \secref{sec:main} by introducing our notion of
covering, and stating our main classification theorems.
Analogues of these theorems for coverings of groupoids were proved in
\cite{rbrown}*{Chapter~9}.
In
\secref{sec:covering action}, we briefly discuss actions of groupoids on
sets, and prove the equivalence between the category of coverings of a
$k$-graph $\Lambda$ and the category of actions of its fundamental
groupoid $\gg\LL$
(Theorem~\ref{equivalence}).
The main theorem of this section is a technical result
(\thmref{groupoid covering})
which implies that the
connected coverings of $\LL$ correspond to transitive actions of $\gg\LL$, and
that the fundamental groups of coverings of $\LL$ can be identified with the
stability groups of the corresponding actions of $\gg\LL$.
In
\secref{sec:analogue}, we state and prove analogues for groupoid actions
of most of our main theorems, and then in \ref{sec:proof main} we prove
the main theorems themselves. Many of them follow from the general
results in the previous section, but when it seemed easier to prove a
result about coverings directly, we did so.

In \secref{sec:skew}, we
construct
universal coverings using skew products. We also show that every
connected covering is a
\emph{relative skew product} (Corollary~\ref{connected
skew}), and prove a version of the Gross-Tucker Theorem which identifies
the $k$-graphs which admit free actions of a group as skew products.

It seems to us that $k$-graphs are likely to be of interest in their own
right, so we have been careful to limit our discussion of $C^*$-algebras to a
final section on the applications of our theory. Our sought-after
generalisation of \cite{dpr}*{Theorem~3.2} is \corref{covering
crossed product}. The main idea in the proof of \corref{covering
crossed product} comes from \cite{kqr:graph}: every group-valued cocycle
$\eta$ on a $k$-graph $\Lambda$ induces a normal and maximal coaction
$\delta_\eta$ of the group on $C^*(\Lambda)$, and every $k$-graph carries a
suitable
cocycle with values in the fundamental group. We also prove a
decomposition theorem which generalises \cite{dpr}*{Corollary~3.6},
prove
that
the $C^*$-algebra of every $k$-graph is nuclear, and
prove
that
the $C^*$-algebra of the skew product by the degree map is always
AF.

Most of this research was done while the second author visited the University
of Newcastle, and he thanks his hosts, particularly Iain Raeburn and David
Pask, for their hospitality.

\section{Main results}
\label{sec:main}

For $k$-graphs and groupoids
we adopt the conventions of \cites{kp:kgraph, pqr:groupoid, rsy}, except
that we do not require them to be countable.
Briefly, a \emph{$k$-graph} is a small category $\LL$ equipped with a
functor $d\:\LL\to\N^k$ satisfying the
\emph{factorization property}: for all $\a\in\LL$ and $n,l\in\N^k$ such
that $d(\a)=n+l$ there exist unique $\b,\g\in\LL$ such that $d(\b)=n$,
$d(\g)=l$, and $\a=\b\g$.
When $d(\a)=n$ we say $\a$ has \emph{degree} $n$.
A \emph{groupoid} is a small category in which every morphism has an
inverse. All groupoids and groups in this paper are discrete, in the
sense that they carry no topology.

If $\C$ is either a $k$-graph or a groupoid,
the \emph{vertices} are the objects,
and $\C^0$ denotes the set of vertices.
For $\a\in\C$, the \emph{source} $s(\a)$ is the domain,
and the \emph{range} $r(\a)$ is the codomain.
For $u,v\in\C^0$ we write $u\C=r\inv(u)$, $\C v=s\inv(v)$, and
$u\C v=u\C \cap
\C v$.
$\C$ is \emph{connected} if the equivalence relation on $\C^0$ generated by
$\{(u,v)\mid u\C v\ne\emptyset\}$ is $\C^0\times\C^0$;
for a groupoid this just means $u\C v\ne\emptyset$ for all $u,v\in \C^0$.
If $\LL$ is a $k$-graph, $u,v\in\LL^0$,
and $n\in\N^k$, we write $\LL^n=d\inv(n)$,
$u\LL^n=u\LL \cap \LL^n$, and
$\LL^n v=\LL v \cap \LL^n$.
A \emph{morphism} between $k$-graphs is a degree-preserving functor.

In general we often write composition of maps as juxtaposition,
especially when we are chasing around commutative diagrams.

\begin{defn}
A \emph{covering} of a $k$-graph $\Lambda$ is a surjective
$k$-graph morphism $p\:\Omega\to \Lambda$ such that
for all $v\in \OO^0$, $p$ maps
$\OO v$ 1-1 onto $\LL p(v)$ and $v\OO$ 1-1 onto $p(v)\LL$.
If $(\OO,p)$ and $(\S,q)$ are coverings of $\LL$,
a \emph{morphism} from $(\OO,p)$ to $(\S,q)$
is a
$k$-graph
morphism $\p\:\OO\to\S$ making the diagram
\[
\xymatrix
@C-5pt
@R-5pt
{
\OO \ar[rr]^-\p \ar[dr]_p
&&{\S} \ar[dl]^q
\\
&\LL
}
\]
commute;
we write $\p\:(\OO,p)\to(\S,q)$.
A covering $p\:\OO\to\LL$ is
\emph{connected} if $\OO$, hence $\LL$, is connected.
\end{defn}

\begin{rem}
If $\LL$ is connected then surjectivity of $p$ is implied by the other
properties.
Also, any functor $\phi\:\Omega\to\Sigma$ making the above diagram
commute automatically preserves degrees, hence is a morphism of coverings.
\end{rem}

Every $k$-graph $\LL$ has a \emph{fundamental groupoid}, which is a groupoid
$\gg\LL$ such that $\gg\LL^0=\LL^0$ together with a canonical functor
$i:\LL\to \gg\LL$ which is the identity on $\LL^0$ and has
the
following universal
property: for every functor $T$ from $\LL$ into a groupoid $\H$
there exists a unique groupoid morphism $T'$ making the diagram
\[
\xymatrix
@C-5pt
@R-5pt
{
\LL \ar[r]^-i \ar[dr]_T
&\gg\LL \ar[d]^{T'}
\\
&\H
}
\]
commute. The assignment
$\LL\mapsto\gg\LL$ is functorial from $k$-graphs to groupoids.
The  \emph{fundamental group}
of $\LL$ at a vertex $x\in\LL^0$ is the isotropy group
\[
\pi(\LL,x):=x\gg\LL x.
\]
(The subscript 1 in the standard notation seems redundant in this context.)
By functoriality of $\LL\mapsto \gg\LL$,
a covering $p\:\OO\to\LL$ induces homomorphisms $p_*\:\pi(\OO,v)\to
\pi(\LL,p(v))$.

\begin{thm}
\label{classify}
Let $(\OO,p)$ and $(\S,q)$ be connected coverings of a $k$-graph $\LL$.
For all $x\in\LL^0$,
the family $\{p_*\pi(\OO,v)\mid p(v)=x\}$ is a conjugacy class
of subgroups of $\pi(\LL,x)$.
For all $x\in\LL^0$, $v\in p\inv(x)$, and $u\in q\inv(x)$,
there is a morphism $(\OO,p)\to(\S,q)$ taking $v$ to $u$ if and only if
$p_*\pi(\OO,v)\subset q_*\pi(\S,u)$.
Consequently, $(\OO,p)\cong(\S,q)$
if and only if
the subgroups $p_*\pi(\OO,v)$
and $q_*\pi(\S,u)$ of $\pi(\LL,x)$ are conjugate
for some, hence every, $x\in\LL^0$, $v\in p\inv(x)$, and $u\in q\inv(x)$.
\end{thm}

For our next result we will need to enlarge our supply of morphisms:

\begin{defn}
If $p\:\OO\to\LL$ and $q\:\OO\to\Gamma$ are coverings,
a \emph{morphism} from $(\OO,p)$ to $(\OO,q)$
is a
$k$-graph
morphism $\p\:\LL\to\Gamma$ making the diagram
\[
\xymatrix
@C-5pt
@R-5pt
{
&\OO \ar[dl]_p \ar[dr]^q
\\
\LL \ar[rr]_-\p
&&\Gamma
}
\]
commute;
we write $\p\:(\OO,p)\to(\OO,q)$.
\end{defn}

From the context it is always clear which type of morphism of coverings we
mean.
There is an obvious notion of morphism which would unify the two kinds we've
introduced, but since we have no use for it we omit it.
In the above definition, observe that since $p$ is surjective,
there is at most one
morphism $\p\:(\OO,p)\to(\OO,q)$.

We will also need to know about quotients by group actions:
let $\aut(\OO,p)$ denote the automorphism group of a connected covering
$p\:\OO\to\LL$.
As we shall show in \secref{sec:proof main},
the quotient map $\OO\to\OO/\aut(\OO,p)$ gives rise to a commuting diagram
\[
\xymatrix
@C-5pt
@R-5pt
{
\OO \ar[d]_p \ar[r]
&\OO/\aut(\OO,p) \ar[dl]
\\
\LL
}
\]
of connected coverings.

\begin{cor}
\label{covering normal}
Let $p\:\OO\to\LL$ be a connected covering, $x\in\LL^0$, and $v\in p\inv(x)$.
Then the following are equivalent:
\begin{enumerate}
\item
the subgroup $p_*\pi(\OO,v)$
of $\pi(\LL,x)$
is normal;

\item
$\aut(\OO,p)$ acts transitively on $p\inv(x)$;

\item
the covering $\OO/\aut(\OO,p)\to\LL$ is an isomorphism;

\item
$(\OO,p)$ is isomorphic to the covering $\OO\to\OO/\aut(\OO,p)$.
\end{enumerate}
\end{cor}

\begin{thm}
\label{covering normalizer}
Let $p\:\OO\to\LL$ be a connected covering, $x\in\LL^0$, and $v\in p\inv(x)$.
Then the normalizer $N(p_*\pi(\OO,v))$ of $p_*\pi(\OO,v)$ in $\pi(\LL,x)$
acts on $(\OO,p)$, and in fact
\[
\aut(\OO,p)\cong N(p_*\pi(\OO,v))/p_*\pi(\OO,v).
\]
\end{thm}

\thmref{classify}
shows how
isomorphism classes of connected coverings of a $k$-graph $\LL$ are
inversely related to conjugacy classes of
subgroups of the fundamental group $\pi(\LL,x)$
(for any choice of vertex $x$).
The identity map on $\LL$ gives a minimal covering,
corresponding to the improper subgroup $\pi(\LL,x)$.
At the opposite extreme:

\begin{defn}
A
covering $p\:\OO\to\LL$ is \emph{universal} if
it is connected and
for every connected covering $q\:\S\to\LL$ there exists a morphism
$(\OO,p)\to(\S,q)$.
\end{defn}

For coverings of small categories, the following result
is \cite{bridson}*{Proposition A.19}.

\begin{thm}
\label{covering universal}
Every connected $k$-graph $\LL$ has a universal covering.
A connected covering
$p\:\OO\to\LL$ is universal if and only if
$p_*\pi(\OO,v)=\{x\}$ for some, hence every,
$x\in\LL^0$ and $v\in p\inv(x)$.
\end{thm}

The following result shows that every subgroup of
$\pi(\LL,x)$ occurs
in the form $p_*\pi(\Omega,v)$ for some
connected covering $(\Omega,p)$:

\begin{thm}
\label{covering quotient}
Let $p\:\OO\to\LL$ be a universal covering, $x\in\LL^0$, $v\in p\inv(x)$, and
$H$ a subgroup of $\pi(\LL,x)$. Let $H$ act on $(\OO,p)$ according to
\thmref{covering normalizer}. Then the associated covering $q\:\OO/H\to\LL$
is connected, and
\[
H=q_*\pi(\OO/H,vH).
\]
Moreover, every connected covering of $\LL$ is isomorphic to one of these
coverings $\OO/H\to\LL$.
\end{thm}

\section{Coverings and actions}
\label{sec:covering action}

\begin{defn}
(cf. \cites{rbrown, higgins})
An \emph{action} of a groupoid $\G$ on a set $V$ is a functor $T$ from $\G$ to
the category of sets such that $V$ is the disjoint union of the sets $T(x)$
for $x\in \G^0$.
Put:
\begin{itemize}
\item
$V_x=T(x)$ for $x\in \G^0$;

\item
$\G*V=\{(a,v)\mid a\in\G,v\in V_{s(a)}\}$;

\item
$av=T(a)(v)$ for $(a,v)\in\G*V$.
\end{itemize}
The \emph{transformation groupoid} is the set $\G*V$ with operations
\begin{align*}
s(a,v)=(s(a),v),
\qquad
r(a,v)=(r(a),av),
\qquad
(a,bv)(b,v)=(ab,v).
\end{align*}
The \emph{stability group}
at $v\in V$ is
\[
S_v:=\{a\in\G\mid av=v\}.
\]
The action $(V,\G)$ is
\begin{itemize}
\item
\emph{transitive} if $V=\G v$ for some (hence every) $v\in V$;

\item
\emph{free} if $S_v=\{x\}$ for all $x\in\G^0$ and $v\in V_x$.
\end{itemize}
\end{defn}

Thus, for each object $x\in\G^0$ we have a set $V_x$, and for each
$a\in x\G y$ we have a bijection $v\mapsto av$ from $V_y\to V_x$.
Since we require the sets
$V_x$ to be pairwise disjoint, we have a bundle $V\to\G^0$, and $\G$ acts
as bijections
among the fibers of this bundle.

\begin{defn}
If $\G$ acts on both $V$ and $U$, a \emph{morphism} from $(V,\G)$
to $(U,\G)$ is a map $\p\:V\to U$ which is \emph{$\G$-equivariant} in the
sense that
\[
\p(av)=a\p(v)
\forall (a,v)\in\G*V.
\]
\end{defn}

\begin{rem}
Thus,
a morphism between actions of $\G$ is just a natural transformation between
the functors.
\end{rem}

We are ready to begin forging the connection between coverings and actions,
which will take the form of
an equivalence between the
categories of coverings of a $k$-graph $\Lambda$ (and their morphisms)
and actions of $\mathcal G(\Lambda)$ (and their morphisms).

\begin{prop}
\label{action to covering}
Let $\LL$ be a $k$-graph, and let its fundamental groupoid $\gg\LL$ act on a
set $V$.
Put
\[
\LL*V=\{(\a,v)\in\LL\times V\mid v\in V_{s(\a)}\}
\]
and
\[
\a v=i(\a)v
\for (\a,v)\in\LL*V.
\]
Then $\LL*V$
becomes a $k$-graph with operations
\begin{align*}
s(\a,v)&=(s(\a),v)
&
r(\a,v)&=(r(\a),\a v)
\\
(\a,\b v)(\b,v)&=(\a\b,v)
&
d(\a,v)&=d(\a),
\end{align*}
and the coordinate projection $p_\LL\:\LL*V\to\LL$ is a covering.

Moreover, the assignments
\begin{equation}
\label{act to cov}
(V,\gg\LL)\mapsto(\LL*V,p_\LL)
\and
\p\mapsto\id_\LL*\p
\end{equation}
give a functor from
actions of $\mathcal G(\Lambda)$ to coverings of $\Lambda$.
\end{prop}

\begin{proof}
Routine computations
(similar to those showing $\G(\LL)*V$ is a groupoid)
verify that $\Lambda*V$ is a category.
Also, $p_\LL$ is clearly a surjective morphism. To see that it has the
covering property, just note that for $(x,v)\in \LL^0*V=(\LL*V)^0$ we have
\begin{align*}
(\LL*V)(x,v)&=\{(\a,v)\mid \a\in\LL x\}
\\
(x,v)(\LL*V)&=\{(\a,u)\mid \a\in x\LL, \a u=v\}.
\end{align*}

The map $d\:\LL*V\to\N^k$ is the composition of
the functors $d\:\LL\to\N^k$ and $p_\LL\:\LL*V\to\LL$, so is a functor.
We verify the factorization property:
let $(\ll,v)\in \LL*V$ and $n,l\in \N^k$
with $d(\ll,v)=n+l$. Then $d(\ll)=n+l$, so there exist $\m\in \LL^n,\n\in
\LL^l$ such that $\ll=\m\n$. Then $(\m,\n v),(\n,v)\in \LL*V$, and since
\[
s(\m,\n v)
=(s(\m),\n v)
=(r(\n),\n v)
=r(\n,v),
\]
we can multiply:
\[
(\m,\n v)(\n,v)=(\m\n,v)=(\ll,v).
\]
Since
\[
d(\m,\n v)=d(\m)=n
\and
d(\n,v)=d(\n)=l,
\]
the factorization has the right degrees. For uniqueness, if also
\[
(\ll,v)=(\a,\b v)(\b,v)\mbox{ with $d(\a,\b v)=n$ and $d(\b,v)=l$,}
\]
then $\ll=\a\b$ with $d(\a)=n$ and $d(\b)=l$, so we must
have $\a=\m$ and $\b=\n$, hence
\[
(\a,\b v)=(\m,\n v)
\and
(\b,v)=(\n,v).
\]
Thus $\LL*V$ is a $k$-graph.
$p_\LL$ preserves degrees by construction, hence is a covering.

For the other part, routine computations show that if $\p\:V\to U$ is a
morphism of actions of $\gg\LL$, then $\id_\LL*\p$ is a morphism of the
corresponding coverings of $\LL$, and that \eqref{act to cov} is functorial.
\end{proof}

In the opposite direction:

\begin{prop}
\label{covering to action}
Let $p\:\OO\to\LL$ be a covering of $k$-graphs.
Then there exists a unique action of
$\gg\LL$ on $\OO^0$ such that
\[
i(\a)v=r(\ll)
\righttext{if}\a=p(\ll)\text{ and }v=s(\ll),
\]
where $i\:\LL\to\gg\LL$ is the canonical functor.

Moreover,
the assignments
\begin{equation}
\label{cov to act}
(\OO,p)\mapsto(\OO^0,\gg\LL)
\and
\p\mapsto\p|\OO^0
\end{equation}
give a functor from
coverings of $\Lambda$ to actions of $\mathcal G(\Lambda)$.
\end{prop}

\begin{proof}
For $\a\in x\LL y$, the range and source maps of $\OO$ take $p\inv(\a)$ 1-1
onto $p\inv(x)$ and $p\inv(y)$, respectively, thus affording a bijection
$T(\a):p\inv(y)\to p\inv(x)$. Routine computations show that the resulting map
$T$ from $\LL$ to the groupoid of bijections among the sets in the family
$\{p\inv(x)\mid x\in\LL^0\}$ is functorial. Since $T$ is a functor from $\LL$
into a groupoid, it factors uniquely through a morphism from $\G(\LL)$ to the
same groupoid, giving the desired action of $\gg\LL$.

For the other part, a routine computation shows that if
$\p\:(\OO,p)\to(\S,q)$
is a morphism of
coverings
of $\LL$,
then $\p(\a v)=\a\p(v)$ for $(\a,v)\in\LL*\OO^0$.
Thus if $T$ and $S$ are the functors giving the actions of $\G(\LL)$ on
$\OO^0$ and $\S^0$, respectively, then
$\p$ gives a natural transformation from $Ti$ to $Si$
(where $i\:\LL\to\G(\LL)$ is the canonical functor),
hence a natural transformation from $T$ to $S$ by universality of $i$.
Therefore $\p|\OO^0$ is a morphism of groupoid actions.
Routine computations show that \eqref{cov to act} is functorial.
\end{proof}

For coverings of groupoids, the following result is \cite{rbrown}*{Section
9.4, Exercise 3}. For coverings of small categories, it is similar to
\cite{bridson}*{Proposition A.23}.

\begin{thm}
\label{equivalence}
The functors described in the preceding two propositions give a category
equivalence
between coverings of $\Lambda$ and actions of $\mathcal G(\Lambda)$.
\end{thm}

\begin{proof}
First, given a covering $(\OO,p)$,
it follows from the definitions that the map
$(p,s)\:\OO\to\LL*\OO^0$
is bijective, and clearly $p=p_\LL(p,s)$.
A computation using the identity
$r(\ll)=p(\ll)s(\ll)$ for $\ll\in\OO$
verifies that $(p,s)$
is functorial. Thus $(\OO,p)\cong(\LL*\OO^0,p_\LL)$.

Next, given an action $(V,\gg\LL)$,
it is obvious that the map from
$V$ to $\LL^0*V=(\LL*V)^0$ taking $v\in V_x$ to
$(x,v)$ is bijective, and it follows
straight
from the definitions that it is
$\LL$-equivariant, hence $\gg\LL$-equivariant. Therefore
$(V,\gg\LL)\cong(\LL^0*V,\gg\LL)$.

If $\p\:(\OO,p)\to(\S,q)$ and $\psi\:(V,\gg\LL)\to(U,\gg\LL)$ are morphisms,
then the diagrams
\[
\xymatrix{
{(\OO,p)}
\ar[r]^-{\cong}
\ar[d]_{\phi}
&{(\LL*\OO^0,p_\LL)}
\ar[d]^{\id_\LL*\phi|\OO^0}
\\
{(\S,q)}
\ar[r]_-{\cong}
&{(\LL*\S^0,p_\LL)}
}
\]
and
\[
\xymatrix{
{(V,\gg\LL)}
\ar[r]^-{\cong}
\ar[d]_{\psi}
&{(\LL^0*V,\gg\LL)}
\ar[d]^{\id_{\LL^0}*\psi}
\\
{(U,\gg\LL)}
\ar[r]_-{\cong}
&{(\LL^0*U,\gg\LL)}
}
\]
commute.
Thus the isomorphisms of the preceding two paragraphs implement a natural
equivalence between the functors of
Propositions~\ref{action to covering} and \ref{covering to action}.
\end{proof}

The above equivalence matches up the automorphism groups:

\begin{cor}
\label{auto iso}
If $p\:\OO\to\LL$ is a covering then the map $\p\mapsto\p|\OO^0$ gives an
isomorphism $\aut(\OO,p)\cong\aut(\OO^0,\gg\LL)$.
\end{cor}

\begin{thm}
\label{groupoid covering}
Let $p\:\OO\to\LL$ be a covering, and let $\G(\LL)$ act on $\OO^0$ as in
\propref{covering to action}.
Then the map
\[
(p_*,s)\:\gg\OO\to\gg\LL*\OO^0
\]
is a groupoid isomorphism.
\end{thm}

Before giving the proof of this theorem, let us
use it to deduce 
the following two results, which are similar to
\cite{rbrown}*{9.4.2}
in the case of coverings of groupoids,
and to \cite{bridson}*{Proposition A.22}
in the case of coverings of small categories.

\begin{cor}
\label{connected}
A covering
$p\:\OO\to\LL$ is
connected if and only if the corresponding groupoid action
$(\OO^0,\gg\LL)$ is transitive.
\end{cor}

\begin{proof}
$\OO$ is connected if and only if $\gg\OO$,
equivalently $\gg\LL*\OO^0$,
is.
For
$(a,v)\in\gg\LL*\OO^0$,
we have
$s(a,v)=(s(a),v)$ and
$r(a,v)=(r(a),av)$.
It follows that $\gg\OO$ is connected if and only if for
all $u,v\in\OO^0$ there exists $a\in\gg\LL$ such that $u=av$, i.e., if and
only if $\gg\LL$ acts transitively on $\OO^0$.
\end{proof}

\begin{cor}
\label{1-1}
Let $p\:\OO\to\LL$ be a covering, $x\in\LL^0$, and $v\in p\inv(x)$.
Then $p_*$ maps the fundamental group $\pi(\OO,v)$ isomorphically onto the
stability group $S_v$
of the corresponding groupoid action $(\OO^0,\gg\LL)$.
\end{cor}

\begin{proof}
The isomorphism of \thmref{groupoid covering} takes the morphism $p_*$ to the
coordinate projection $\gg\LL*\OO^0\to\gg\LL$.
It follows that $p_*$ maps $\gg\OO v$ 1-1 onto $\gg\LL x$.
Thus $p_*$ maps $\pi(\OO,v)$
isomorphically onto \emph{some} subgroup of $\pi(\LL,x)$.
For $c\in\gg\OO v$ and $a=p_*(c)$
we have
$av=r(c)$,
so $c\in\pi(\OO,v)$ if and only if $a\in S_v$. The result follows.
\end{proof}

In \cite{bridson}*{Proposition A.17} (for coverings of small categories),
the above injectivity of $p_*$ is asserted to follow ``directly from the
definition of a covering'', but to us it doesn't seem so immediate.

\begin{proof}[Proof of \thmref{groupoid covering}]
Our strategy is to present
the fundamental group\-oids of $\Omega$ and $\Lambda$
as path categories of augmented graphs modulo cancellation
relations and commuting squares, and match up the kernels.
More precisely, we will build a commutative diagram
\[
\xymatrix{
\ff{F^+} \ar[r]^-R_-{\text{onto}} \ar[d]_{(q_*,s)}^\cong
&\gg\OO \ar[d]^{(p_*,s)}
\\
\ff{E^+}*\OO^0 \ar[r]_{Q*\id}^-{\text{onto}}
&\gg\LL*\OO^0
}
\]
of functors, where the left-hand vertical isomorphism takes the equivalence
relation determined by the top horizontal functor $R$ onto the equivalence
relation determined by the bottom horizontal functor $Q*\id$. This will
suffice to show that $(p_*,s)$ is an isomorphism, since the horizontal
functors are surjective.

Let $E$ be the \emph{1-skeleton} of $\LL$,
that is, the graph whose vertices coincide with those of $\LL$
and whose edges $E^1$ comprise all elements of $\LL$ whose degree is a
standard basis vector in $\N^k$.
We need to recall a few things from \cite{pqr:groupoid}.
A \emph{diagram of type $E$} in a category $\C$ is a
map $D\:E\to\C$ which is a
morphism from $E$ to the underlying graph of $\C$.
There is a (small) \emph{path category} $\ff E$
and a \emph{canonical diagram} $\DD\:E\to\ff E$
with the universal property that for every diagram $D\:E\to\C$
there is a unique functor $T$ making the diagram
\[
\xymatrix
@C-5pt
@C-5pt
{
E \ar[r]^-\DD \ar[dr]_D
&\ff E \ar[d]^T
\\
&\C
}
\]
commute.
The assignment $E\mapsto \ff E$ is functorial from graphs to small
categories.
A \emph{relation} for $E$ is a pair $(\a,\b)$ of paths in $\ff E$ with
$s(\a)=s(\b)$ and $r(\a)=r(\b)$. If $K$ is a set of relations for $E$, a
diagram $D$ of type $E$ \emph{satisfies $K$} if $D(\a)=D(\b)$ for all
$(\a,\b)\in K$.
Let $S_\LL$ denote the set of all \emph{commuting squares} for $\LL$, i.e.,
relations for $E$ of the form $(ef,gh)$,
where
$e$ and $f$ are composable edges in $E$ with orthogonal degrees and
$g$ and $h$ are the unique edges
such that
$d(g)=d(f)$, $d(h)=d(e)$ and $ef=gh$.

The \emph{augmented graph}
$E^+=E\cup E\inv$,
where $E\inv$ denotes the inverse edges, can be used to give a presentation of
the fundamental groupoid; more precisely,
letting
$C_E$ be the set $\{(e\inv e,s(e))\mid e\in E^1\cup E\inv\}$ of
\emph{cancellation relations} for $E$,
there is
a surjective functor, which we denote for this proof by $Q$,
making the diagram
\[
\xymatrix{
E^+ \ar[d]_\DD
&E \ar@{_(->}[l] \ar[d]^{i|E}
\\
\ff{E^+} \ar[r]^{\text{onto}}_-Q
&\gg\LL
}
\]
commute,
such that the associated equivalence relation on $\ff{E^+}$
is generated by $C_E\cup S_\LL$,
where $\DD\:E^+\to\ff{E^+}$ is the canonical diagram and
$i\:\LL\to\gg\LL$ is the canonical functor.
In particular, the diagram $Q\DD\:E^+\to\gg\LL$ satisfies $C_E$, i.e.,
\[
Q\DD(e\inv)=Q\DD(e)\inv
\forall e\in E^1.
\]
Let $R\:\ff{F^+}\to\gg\OO$ be the corresponding
surjective functor for the 1-skeleton $F$
of the covering $k$-graph $\OO$.

Consider the diagram
\[
\xymatrix
@C-5pt
@R-5pt
{
\ff{F^+} \ar[rrrr]^-R \ar[ddd]_{q_*}
&&&&\gg\OO \ar[ddd]^{p_*}
\\
&F^+ \ar[ul]^\DD \ar[d]_q
&F \ar@{_(->}[l] \ar@{^(->}[r] \ar[d]^{p|F}
&\OO \ar[ur]_i \ar[d]^p
\\
&E^+ \ar[dl]_\DD
&E \ar@{_(->}[l] \ar@{^(->}[r]
&\LL \ar[dr]^i
\\
\ff{E^+} \ar[rrrr]_-Q
&&&&\gg\LL
}
\]
The $\DD$'s are canonical diagrams and the $i$'s are canonical functors.
The right-hand quadrilateral commutes by functoriality of
$\LL\mapsto\G(\LL)$.
The restriction $p|F$ is a graph morphism,
and takes $F$ onto $E$ by
definition of covering and 1-skeleton.
The graph morphism $q\:F^+\to E^+$ is
the extension of $p|F$ defined by
\[
q(f\inv)=p(f)\inv.
\]
The inside squares commute by definition of $p|F$ and $q$.
The left-hand quadrilateral commutes by functoriality of
$E\mapsto\ff E$.
The top and bottom 5-sided diagrams
commute by construction of $R$ and $Q$.
Thus
\[
p_*R\DD|F=Qq_*\DD|F.
\]
Since $Q\DD$ satisfies $C_E$, $Q\DD q=Qq_*\DD$ satisfies $C_F$;
since $R\DD$ also
satisfies $C_F$ and $p_*$ is a groupoid morphism, it follows that
$p_*R\DD=Qq_*\DD$, hence $p_*R=Qq_*$ by universality of $\DD$.

Thus the diagram
\[
\xymatrix{
\ff{F^+} \ar[r]^-R \ar[d]_{q_*}
&\gg\OO \ar[d]^{p_*}
\\
\ff{E^+} \ar[r]_-Q
&\gg\LL
}
\]
commutes.
It is easy to see that the slightly enlarged diagram
\[
\xymatrix{
\ff{F^+} \ar[r]^-R \ar[d]_{(q_*,s)}
&\gg\OO \ar[d]^{(p_*,s)}
\\
\ff{E^+}*\OO^0 \ar[r]_-{Q*\id}
&\gg\LL*\OO^0
}
\]
also commutes --- the sources just come along for the ride.
This is the
diagram indicated at the beginning of the proof.
$q_*\:\ff{F^+}\to\ff{E^+}$ is a 1-graph covering, and
the map $(q_*,s)$ is the associated isomorphism
from the proof of \thmref{equivalence}.

The horizontal functors $R$ and $Q*\id$ are surjective,
so it remains to prove that
$(q_*,s)$ takes the equivalence relation determined by $R$,
which is generated by
the cancellation relations $C_E$ and the commuting squares $S_\Omega$,
onto the equivalence
relation determined by $Q*\id$, which is generated
by the pairs of the form 
$((\a,v),(\b,v))$, where $(\a,\b)\in C_E\cup S_\LL$.

Let $(f^{-1}f,s(f))\in C_F$, and put $e=q_*(f)\in E^+$.
We have
\[
q_*(f\inv f)=q_*(f)\inv q_*(f)=e\inv e,
\]
and $s(f\inv f)=s(f)$, so
\[
(q_*,s)(f\inv f)=(e\inv e,s(f)).
\]
On the other hand,
$q_*(s(f))=s(e)$,
so
\[
(q_*,s)(s(f))=(s(e),s(f)).
\]
Thus
\begin{equation}
\label{relation 1}
(q_*,s)(f^{-1}f,s(f))=\bigl((e\inv e,s(f)),(s(e),s(f))\bigr);
\end{equation}
note that $(e^{-1}e,s(e))$ is a typical element of $C_E$.

Now let 
$(ab,cd)\in S_\Omega$, so that $a,b,c,d\in F^1$, 
the common degree of $a$ and $d$ is orthogonal to the
common degree of $b$ and $c$, and $ab=cd$ in the $k$-graph $\OO$.
Put
\[
e=q(a),
\quad
f=q(b),
\quad
g=q(c),
\and
h=q(d).
\]
Then
\[
q_*(ab)=ef \and q_*(cd)=gh,
\]
and $ef=gh$ in the $k$-graph $\LL$ because the diagram
\[
\xymatrix
@C-5pt
@R-5pt
{
F \ar@{^(->}[r] \ar[d]_{p|F}
&\OO \ar[d]^p
\\
E \ar@{^(->}[r]
&\LL
}
\]
commutes.
On the other hand,
$s(ab)=s(cd)$,
so
\begin{equation}
\label{relation 2}
(q_*,s)(ab,cd)=\bigl((ef,s(ab)),(gh,s(ab))\bigr);
\end{equation}
since $q:F\to E$ is a covering, $(ef,gh)$ is a typical element of $S_\Lambda$.

Together, Equations~\eqref{relation 1} and \eqref{relation 2}
show that $(q_*,s)$ takes $C_F\cup
S_\Omega$ onto the set
$\{((\a,v),(\b,v))\mid (\a,\b)\in C_E\cup S_\LL\}$,
which suffices.
\end{proof}

It follows from \thmref{groupoid covering} that if
$p\:\OO\to\LL$ is a $k$-graph covering,
then $p_*\:\gg\OO\to\gg\LL$ is a \emph{groupoid covering}
in the sense of \cites{rbrown, higgins}.

\section{Classification of transitive groupoid actions}
\label{sec:analogue}

By the results of the preceding section, to classify connected coverings we
only need the well-known classification of transitive groupoid actions. The
results we state in this section are elementary and we claim no originality.
We supply the proofs for the convenience of the reader.

\begin{prop}
\label{action conjugacy}
Let $(V,\G)$ be a transitive groupoid action and $x\in\G^0$.
Then the family $\{S_v\mid v\in V_x\}$ is a conjugacy class
of subgroups of $ x\G x $.
\end{prop}

\begin{proof}
Just note that
$S_{av}=aS_va\inv$ for
$(a,v)\in\G*V$.
\end{proof}

\begin{prop}
\label{action morphism}
Let a groupoid $\G$ act transitively on both $V$ and $U$,
and let $x\in\G^0$, $v\in V_x$, and $u\in U_x$.
Then there is a morphism $(V,\G)\to(U,\G)$ taking $v$ to $u$ if and only if
$S_v\subset S_u$.
\end{prop}

\begin{proof}
If
$\p\:V\to U$ is equivariant and $\p(v)=u$, then
\[
\p(av)=a\p(v)=au
\forall
a\in\G x,
\]
so $S_v\subset S_u$.
Conversely, assume $S_v\subset S_u$, and define $\p\:V\to U$ by
\[
\p(av)=au
\for
a\in\G x.
\]
This is well-defined because if $a,b\in\G x$ and $av=bv$, then
\[
b\inv a\in S_v\subset S_u,
\]
so $au=bu$. Clearly $\p$ is equivariant.
\end{proof}

\begin{prop}
\label{action normalizer}
Let $(V,\G)$ be a transitive groupoid action, $x\in\G^0$, and $v\in V_x$.
Then the normalizer $N(S_v)$ of $S_v$ in $ x\G x $ acts on
the right of
the action $(V,\G)$ by automorphisms, and in fact
\[
\aut(V,\G)\cong N(S_v)/S_v.
\]
\end{prop}

\begin{proof}
The computations in the proof of \propref{action morphism} show that every
automorphism of $(V,\G)$ is of the form
$av\mapsto acv$,
where $c\in x\G x $ satisfies
\[
S_v=S_{cv}=cS_vc\inv,
\]
i.e., $c\in N(S_v)$, and conversely every such $c$ gives rise to an
automorphism of $(V,\G)$ in this manner. Define
$(av)c=acv$.
Then $N(S_v)$ acts on the right, since for $c,d\in N(S_v)$ we have
\[
\bigl((av)c\bigr)d=(acv)d=acdv=(av)cd.
\]
Clearly $(av)c=av$ if and only if $c\in S_v$. Thus $N(S_v)/S_v$ acts freely on
the right of $(V,\G)$. The result follows.
\end{proof}

Our next result is the groupoid-action analogue of \thmref{covering
universal}. However, rather than merely asserting the existence of a certain
kind of action of a groupoid, we give more detail, because this will be useful
when we apply it to the analogue for coverings.
First, we need:

\begin{defn}
A \emph{cocycle} on a groupoid $\G$ is a functor
$\eta\:\G\to G$ where $G$ is a group.
The \emph{cocycle action} of $\G$ on the Cartesian product $\G^0\times G$
is given
by
\[
a(s(a),g)=(r(a),\e(a)g).
\]
We write $\G^0\times_\e G$ to indicate $\G^0\times G$ equipped
with the cocycle action.
\end{defn}

\begin{prop}
\label{action universal}
Let $\G$ be a connected groupoid and $x\in\G^0$. There is a cocycle
$\e\:\G\to x\G x $ such that the associated cocycle action
$(\G^0\times_\e x\G x ,\G)$ is free and transitive.
\end{prop}

\begin{proof}
For each $y\in\G^0$ pick $t_y\in y\G x$, with $t_x=x$.
Then $\e(a)=t_{r(a)}\inv at_{s(a)}$ defines a surjective cocycle
$\G\to x\G x$ which is the identity map on $x\G x$.

For $a\in y\G z$ and $g\in x\G x $ we have
$a(z,g)=(y,\e(a)g)$.
The action is transitive, because if $y\in\G^0$ and $g\in x\G x $ we have
$t_yg(x,x)=(y,g)$.
By transitivity,
to show that the action is free
it suffices to check the stability group at
$(x,x)$: for $a\in y\G x$ we have $a(x,x)=(y,\e(a))$, so if $a\in
S_{(x,x)}$ then $y=x$ and $\e(a)=x$, so $a\in x\G x $, hence $\e(a)=a$, thus
$a=x$.
\end{proof}

We next give a groupoid-action analogue of \thmref{covering
quotient}.
First note that if a group $G$ acts on (the
right of) a groupoid action $(V,\G)$ by automorphisms,
then $\G$ acts on the quotient set $V/G$ by
$a(vG)=(av)G$.

\begin{prop}
\label{action quotient}
Let $(V,\G)$ be a free transitive groupoid action, $x\in\G^0$, $v\in V_x$, and
$H$ a subgroup of $ x\G x $. Let $H$ act on $(V,\G)$ according to
\propref{action normalizer}. Then the associated action $(V/H,\G)$
is transitive, and
$H=S_{vH}$.
\end{prop}

\begin{proof}
The action of $\G$ on $V/H$ is transitive since it is a quotient of the
transitive action on $V$.
We have $V=\G v$, $H$ acts on $V$ by
$(av)h=ahv$,
and $\G$ acts on $V/H$ by
$a(vH)=avH$.
Thus $a(vH)=vH$ if and only if $av\in Hv$, equivalently $a\in H$ by freeness.
\end{proof}

\section{Proofs of main results}
\label{sec:proof main}

\begin{proof}[Proof of \thmref{classify}]
The corresponding groupoid actions
$(\OO^0,\gg\LL)$ and $(\S^0,\gg\LL)$
are transitive, by \corref{connected}.
We have
$\pi(\LL,x)= x\gg\LL x $, and for $v\in p\inv(x)$ we have $p_*\pi(\OO,v)=S_v$,
so the first statement follows from \propref{action conjugacy}.

By \thmref{equivalence}
a morphism
$\p\:(\OO,p)\to(\S,q)$ with $\p(v)=u$
corresponds to a morphism
$\psi\:(\OO^0,\gg\LL)\to(\S^0,\gg\LL)$
with $\psi(v)=u$, and we have $p_*\pi(\OO,v)=S_v$ and $q_*\pi(\S,u)=S_u$, so
the second statement follows from \propref{action morphism}.

The last statement now follows quickly from the above.
\end{proof}

\begin{prop}
Let $G$ be a group acting freely by automorphisms on the right of a $k$-graph
$\OO$. Then the quotient set $\OO/G$ becomes a $k$-graph with operations
\begin{align*}
s(\ll G)&=s(\ll)G
&
r(\ll G)&=r(\ll)G
\\
(\ll G)(\m G)&=(\ll\m)G
&
d(\ll G)&=d(\ll),
\end{align*}
and the quotient map $\OO\to\OO/G$ is a covering.
\end{prop}

\begin{proof}
More precisely, the composition is defined as follows:
if $s(\ll)G=r(\m)G$, then it follows from
freeness of the action that the set
\[
\{\a\b\mid \a\in \ll G,\b\in \m G,s(\a)=r(\b)\}
\]
comprises a single orbit,
which coincides with $(\ll\m)G$ if we
adjust $\ll,\m$ within their respective orbits so
that
they are composable
(and freeness is needed to show that the composition of
orbits is well-defined).
Routine computations show that $\OO/G$ is a
category, and the quotient map is then a surjective functor.
Because $G$ acts by automorphisms, all elements of any orbit $\a G$ have the
same degree, which we define to be $d(\a G)$. This gives a functor
$d\:\OO/G\to\N^k$, and
by construction
the quotient map intertwines the two ``$d$''s.

We verify the factorization property: 
let $\ll G\in\OO/G$ and $n,l\in\N^k$ with $d(\ll
G)=n+l$. Then $d(\ll)=n+l$, so there exist
unique
$\m\in\OO^n$ and $\n\in\OO^l$ such
that $\ll=\m\n$, and then
\[
d(\m G)=n,
\quad
d(\n G)=l,
\and
\ll G=\m G\n G.
\]
For the uniqueness of $\ll G$ and $\m G$, suppose
\[
d(\a G)=n,
\quad
d(\b G)=l,
\and
\ll G=\a G\b G.
\]
Then $d(\a)=n$ and $d(\b)=l$, and we can
adjust $\b$
in the $G$-orbit
so that
$s(\a)=r(\b)$.
Then
$\a G\b G=(\a\b)G$,
so there exists $g\in G$ such that
\[
(\a g)(\b g)=(\a\b)g=\ll.
\]
But
$d(\a g)=n$
and
$d(\b g)=l$,
so we must have
$\a g=\m$ and $\b g=\n$, hence
$\a G=\m G$ and $\b G=\n G$.

Thus $\OO/G$ is a $k$-graph, and the quotient map is a $k$-graph morphism.
For $v\in\OO^0$,
the quotient map takes $\OO v$ onto $(\OO/G)(vG)$ by construction;
we must show that it is injective on this set.
Let $\ll,\m\in\OO v$ such that $\ll G=\m G$. Then there exists
$g\in G$ such that $\ll=\m g$. Thus
\[
v=s(\ll)=s(\m g)=s(\m)g=vg.
\]
Since $G$ acts freely, we must have
$g=e$, hence $\ll=\m$. Similarly for $v\OO$ and $(vG)(\OO/G)$.
\end{proof}

\begin{prop}
\label{associated}
Let $G$ be a group acting freely by covering automorphisms on the right of a
$k$-graph covering
$p\:\OO\to\LL$. Then the map
$\ll G\mapsto p(\ll)\:\OO/G\to\LL$
is a covering.
\end{prop}

\begin{proof}
We could deduce this from a corresponding groupoid-action result, but it's
faster to prove this one directly. We certainly have a commuting diagram
\[
\xymatrix
@C-10pt
@R-10pt
{
\OO \ar[r] \ar[d]_p
&\OO/G \ar[dl]
\\
\LL
}
\]
of surjective functors,
where $\OO\to\OO/G$ is the covering
from the preceding proposition.
It is easy to
verify that, whenever we have such a commuting diagram of surjective functors,
if two of the maps are coverings, then so is the third.
\end{proof}

For the proof of our next main result,
we wish to apply the above to
the automorphism group of a
connected covering. For this we need to know that this group acts
freely:

\begin{prop}
\label{covering auto free}
Every automorphism of a connected $k$-graph covering acts freely.
\end{prop}

\begin{proof}
Since the covering is connected, the corresponding groupoid action is
transitive, and it is straightforward to verify that every automorphism of a
transitive groupoid action acts freely.
\end{proof}

\begin{proof}[Proof of \corref{covering normal}]
The equivalence (i) $\iff$ (ii) follows quickly from
\thmref{classify}.

The quotient $k$-graph $\OO/\aut(\OO,p)$ is connected since $\OO$ is, hence
the covering $q\:\OO/\aut(\OO,p)\to\LL$ is connected.
Thus this covering is an isomorphism if and only if $q\inv(x)=\{vG\}$.
Since the set $q\inv(x)$ coincides with the set of
$\aut(\OO,p)$-orbits of elements of $p\inv(x)$, we have
(ii) $\iff$ (iii).

Finally, for (iii) $\iff$ (iv), just
note that the covering $\OO/\aut(\OO,p)\to\LL$ is the unique
morphism from the covering $\OO\to\OO/\aut(\OO,p)$ to the given covering
$\OO\to\LL$.
\end{proof}

\begin{proof}[Proof of \thmref{covering normalizer}]
Passing to the corresponding groupoid action $(\OO^0,\gg\LL)$, the result
follows from \propref{action normalizer}.
\end{proof}

Later we will need a precise description of the action of
$N(p_*\pi(\OO,v))$ on $(\OO,p)$ corresponding to the action of $N(S_v)$ on
$(\OO^0,\gg\LL)$, and we record this here: let $\ll\in\OO$ and $c\in
N(p_*\pi(\OO,v))$. Since the covering $(\OO,p)$ is connected, by
\corref{connected} $\gg\LL$ acts transitively on $\OO^0$, so there exists
$a\in\gg\LL$ such that $s(\ll)=av$. Then $\ll c$ is the unique element of
$\OO$ such that
\[
p(\ll c)=p(\ll)
\and
s(\ll c)=acv.
\]
This is well-defined because $c$ normalizes $p_*\pi(\OO,v)$.

\begin{proof}[Proof of \thmref{covering universal}]
By \cite{pqr:groupoid}*{Proposition 5.9}, the fundamental group\-oid $\gg\LL$
is connected since $\LL$ is. \propref{action universal} gives a certain free
and transitive action $(V,\gg\LL)$. Let $p\:\OO\to\LL$ be the corresponding
covering, which is connected by \corref{connected}, and let $x\in\LL^0$ and
$v\in p\inv(x)$. Then
\[
p_*\pi(\OO,v)=S_v=\{x\},
\]
so the covering $(\OO,p)$ is universal by \thmref{classify}.

Moreover, again by \thmref{classify}, if $(\S,q)$ is any universal
covering of $\LL$, then because there is a morphism $(\S,q)\to(\OO,p)$, we
must have $q_*\pi(\S,u)=\{x\}$ for all $u\in q\inv(x)$.
\end{proof}

It will be useful to record the following alternative characterization of
universal coverings. But first:

\begin{defn}
A \emph{$k$-tree} is a connected $k$-graph $\OO$ with
$\pi(\OO,v)=\{v\}$ for some, hence every, vertex $v\in\OO^0$.
\end{defn}

\begin{cor}
\label{free}
If
$p\:\OO\to\LL$ is a connected covering, then the following are equivalent:
\begin{enumerate}
\item
the covering $(\OO,p)$ is universal;

\item
the
corresponding groupoid action $(\OO^0,\gg\LL)$ is free;

\item
$\OO$ is a $k$-tree.
\end{enumerate}
\end{cor}

\begin{proof}
This follows from
\thmref{covering universal}
since $p_*|\pi(\OO,v)$ is injective and $S_v=p_*\pi(\OO,v)$.
\end{proof}

\begin{rem}
A $1$-tree is the path category
(modulo conventions regarding composition)
of a graph which is a tree in the
usual sense.
In a 1-tree, between any 2 vertices there is at most 1
undirected path, hence certainly at most 1 directed path;
this
does not generally hold in $k$-trees, as illustrated by one of our basic
examples \cite{pqr:groupoid}*{Example 7.2} of a $k$-graph $\LL$ which does not
embed faithfully in its fundamental groupoid. This was a 2-graph with 4
vertices, 4 horizontal edges, and 6 vertical edges. All multiple edges between
vertices collapsed under the canonical functor $i$, making the fundamental
groupoid an equivalence relation on 4 objects, thus the fundamental groups
were all trivial. Hence this is an example of a 2-tree with multiple morphisms
with the same source and range.
This is unfortunate, because it means that in practice we have no effective
algorithm for determining whether a given $k$-graph is a $k$-tree, short of
computing the fundamental group.
\end{rem}

\begin{proof}[Proof of \thmref{covering quotient}]
First of all, since the covering is universal,
we have $p_*\pi(\OO,v)=\{x\}$.
Thus
the action of $H$ on $(\OO,p)$
guaranteed by \thmref{covering normalizer} is free,
hence
by \propref{associated}
we really do have a
covering $q\:\OO/H\to\LL$.
Moreover, this covering is connected since it is a
quotient of the connected covering $(\OO,p)$.
By \corref{free} the corresponding groupoid action $(\OO^0,\gg\LL)$ is free.
We have $x\gg\LL x=\pi(\LL,x)$, so $H$ acts on $(\OO^0,\gg\LL)$ according to
\propref{action normalizer}.
By \propref{action quotient} we have $H=S_{vH}$.
Since $q_*\pi(\OO/H,vH)=S_{vH}$ by \corref{1-1}, we have shown the first part
of the theorem.
The other part now follows immediately from \thmref{classify}.
\end{proof}

\section{Skew products}
\label{sec:skew}

Our statement of \thmref{covering universal} merely asserted the existence of
universal coverings, and the proof merely showed how this existence followed
from the analogous
\propref{action universal} for groupoid actions. However,
\propref{action universal} gave a
specific
construction of the desired groupoid action
using cocycles.
The analogue for $k$-graph coverings
is a skew product:

\begin{defn}
A \emph{cocycle} on a $k$-graph $\LL$ is  a functor from $\LL$ to a group.
\end{defn}

\begin{obs}
Since a cocycle $\e\:\LL\to G$ is a functor into a group\(oid\),
there is a unique cocycle $\k$ making the diagram
\[
\xymatrix
@C-5pt
@R-5pt
{
\LL \ar[r]^-\e \ar[d]_i
&G
\\
\gg\LL \ar[ur]_\k
}
\]
commute.
In fact, this gives a 1-1 correspondence between cocycles on $\LL$ and on
$\gg\LL$.
\end{obs}

\begin{prop}
\label{skew product covering}
Let $\e\:\LL\to G$ be a $k$-graph cocycle. Then the Cartesian product
$\LL\times G$ becomes a $k$-graph with operations
\begin{align*}
s(\a,g)&=(s(\a),g)
&
r(\a,g)&=(r(\a),\e(\a)g)
\\
(\a,\e(\b)g)(\b,g)&=(\a\b,g)
&
d(\a,g)&=d(\a),
\end{align*}
and the coordinate projection $\LL\times G\to\LL$ is a covering.
\end{prop}

\begin{proof}
Let $\e'\:\G(\LL)\to G$ be the corresponding groupoid cocycle.
The associated cocycle action of $\G(\LL)$
is on the set $\LL^0\times G$, and the
covering $k$-graph
corresponding to this groupoid action
is $\LL*(\LL^0\times G)$. The map
\[
\bigl(\a,(x,g)\bigr)\mapsto(\a,g)\:\LL*(\LL^0\times G)\to\LL\times G
\]
is bijective, transforms the $k$-graph operations on $\LL*(\LL^0\times G)$
into the operations on $\LL\times G$ indicated in the proposition, and
transforms the corresponding
covering $\LL*(\LL^0\times G)\to\LL$ into the coordinate projection $\LL\times
G\to\LL$.
\end{proof}

The following definition is a variation of \cite{kp:kgraph}*{Definition 5.1}:

\begin{defn}
The \emph{skew product} $k$-graph associated to a cocycle $\e\:\LL\to G$ is
$\LL\times G$ with the operations from \propref{skew product covering}.
We write $\LL\times_\e G$ to indicate this $k$-graph.
The \emph{skew-product covering} is the coordinate projection
$p_\LL\:\LL\times_\e G\to \LL$.
\end{defn}

Note that in the proof of the above proposition, the skew-product covering
$\LL\times_\e G\to\LL$ was not exactly the same as the covering corresponding
to the cocycle action $(\G(\LL)^0\times_\eta G,\G(\LL))$, rather these
coverings were merely
isomorphic --- nevertheless for convenience we regard the skew-product
covering as \emph{corresponding} to the cocycle action, thus committing a mild
abuse.

We can now apply this to construct universal coverings:

\begin{cor}
\label{skew universal}
Let $\LL$ be a connected $k$-graph and $x\in\LL^0$.
Then there is a cocycle
$\e\:\LL\to\pi(\LL,x)$
such that
the skew-product covering
$\LL\times_\e\pi(\LL,x)\to\LL$ is universal.
\end{cor}

\begin{proof}
In
the proof of \propref{action universal} we 
constructed a groupoid cocycle $\gg\LL\to\pi(\LL,x)$;
let
$\e\:\LL\to\pi(\LL,x)$
be the associated $k$-graph cocycle.
The skew-product covering corresponds to
the cocycle action of $\gg\LL$,
and
the proof of \thmref{covering universal} showed
that this corresponding covering is universal.
\end{proof}

\begin{prop}
\label{skew action}
Let $\e\:\LL\to G$ be a $k$-graph cocycle. Then $G$ acts
freely
on the skew-product
covering $\LL\times_\e G\to\LL$ via
\[
(\ll,g)h=(\ll,gh)
\for \ll\in\LL,g,h\in G.
\]
\end{prop}

\begin{proof}
This can be checked directly without pain, but it is even easier to note that
$G$ acts on the corresponding groupoid action $(\LL^0\times G,\gg\LL)$ by
\[
(x,g)h=(x,gh)
\for x\in\LL^0,g,h\in G,
\]
and then apply \corref{auto iso}. It is obvious that the action is free.
\end{proof}

For a connected $k$-graph $\LL$, \corref{skew universal} gives a specific
construction of a universal covering $\LL\times_\e\pi(\LL,x)\to\LL$, so
\thmref{covering normalizer} gives an action of $\pi(\LL,x)$ on this covering.
In the following result we verify that this action coincides with the one
guaranteed by \propref{skew action}:

\begin{prop}
\label{skew normalizer}
Let $\LL$ be a connected $k$-graph, and let $\LL\times_\e\pi(\LL,x)\to\LL$ be
the universal covering as in \corref{skew universal}. Then the action of
$\pi(\LL,x)$ on the skew product covering $\LL\times_\e\pi(\LL,x)\to \LL$ given by
\thmref{covering normalizer} agrees with the action given by \propref{skew
action}.
\end{prop}

\begin{proof}
The corresponding action of $\gg\LL$ is on
$(\LL\times_\e \pi(\LL,x))^0=\LL^0\times \pi(\LL,x)$.
Let's see how
the proof of \thmref{covering normalizer}
tells us
$\pi(\LL,x)$ acts on $\LL\times_\e \pi(\LL,x)\to \LL$:
denoting the skew product covering map $\LL\times_\e \pi(\LL,x)\to\LL$
by $p_\LL$,
we must start by choosing a vertex $v\in p_\LL\inv(x)$. Then the action of an
element $h\in \pi(\LL,x)$ on an element $(\ll,g)\in\LL\times_\e \pi(\LL,x)$ is
computed as follows: find $b\in\gg\LL x$ such that
$bv=s(\ll,g)$,
and then 
$(\ll,g)h$ is the unique element of $\LL\times_\e \pi(\LL,x)$ such that both
\[
p_\LL\bigl((\ll,g)h\bigr)=p_\LL(\ll,g)=\ll
\and s\bigl((\ll,g)h\bigr)=bhv.
\]
For $v$ we choose $(x,x)$.
By definition of the cocycle action,
for any $b\in\gg\LL x$ we have $b(x,x)=\bigl(r(b),\e(b)\bigr)$.
Given $(\ll,g)$, put $y=s(\ll)$. Then
$s(\ll,g)=(y,g)$, so we want $b\in\gg\LL x$ such that
\[
\bigl(r(b),\e(b)\bigr)=(y,g).
\]
The cocycle $\e$ constructed in the proof of \propref{action universal} takes
$y\gg\LL x$ onto $\pi(\LL,x)$, so such an element $b$ exists.
Then for such a $b$ we have
\[
bh(x,x)
=\bigl(r(bh),\e(bh)\bigr)
=\bigl(r(b),\e(b)\e(h)\bigr)
=(y,gh).
\]
Therefore
$(\ll,g)h=(\ll,gh)$,
as desired, since
\[
p_\LL(\ll,gh)=\ll
\and
s(\ll,gh)=(y,gh).
\qedhere
\]
\end{proof}

\thmref{covering quotient} concerns an action of a subgroup $H$ of the
fundamental group $\pi(\LL,x)$ on a universal covering of $\LL$; we want to
see how this looks when $\LL$ is a skew product $\LL\times_\e\pi(\LL,x)$ as in
\corref{skew universal}. It is cleaner to do it in the abstract: let
$\e\:\LL\to G$ be a cocycle and $H$ a subgroup of $G$, and let $H$ act on the
skew product $\LL\times_\e G$ according to \propref{skew action}. Since this
action is free, we can form the associated covering
$(\LL\times_\e G)/H\to\LL$.
The map
\[
(\ll,g)H\mapsto(\ll,gH)
\:
(\LL\times_\e G)/H\to \LL\times (G/H)
\]
is bijective, transforms the $k$-graph operations on
the quotient
$(\LL\times_\e G)/H$ into
\begin{align*}
s(\ll,gH)&=(s(\ll),gH)
&
r(\ll,gH)&=(r(\ll),\e(\ll)gH)
\\
(\ll,\e(\m)gH)(\m,gH)&=(\ll\m,gH)
&
d(\ll,gH)&=d(\ll),
\end{align*}
and transforms the covering $(\LL\times_\e G)/H\to\LL$ into the coordinate
projection $\LL\times G/H\to\LL$.

\begin{defn}
If $\e\:\LL\to G$ is a $k$-graph cocycle and $H$ is a subgroup of $G$,
the \emph{relative skew product} $k$-graph, denoted
$\LL\times_\e G/H$, is 
the Cartesian product $\LL\times G/H$ with the above
operations, and
the \emph{relative skew product covering} is the coordinate
projection $\LL\times_\e G/H\to\LL$.
\end{defn}

We should point out that this concept is not new: a version of relative skew
products for graphs appears in, for example, \cites{gt,
dpr}. While we did not need relative skew products for the general theory of
coverings --- for us they arose as just a particular case of quotients of skew
products --- they will be important for us in our application to
$C^*$-coactions.

Let's formalize the above discussion:

\begin{prop}
\label{orbit skew}
If $\e\:\LL\to G$ is a $k$-graph cocycle and $H$ a subgroup of $G$, then
the associated covering $(\LL\times_\e G)/H\to\LL$ is
isomorphic to the relative skew product covering $\LL\times_\e G/H\to\LL$ via
the map
\[
(a,g)H\mapsto (a,gH).
\]
\end{prop}

The value of the above definition is that it captures all connected coverings,
as we show in the following result, a graph version of which appeared
in \cite{dpr}*{Proposition 2.2}:

\begin{cor}
\label{connected skew}
Every connected covering is isomorphic to a relative skew-product covering.
\end{cor}

\begin{proof}
Let $p\:\OO\to\LL$ be a connected covering, $x\in\LL^0$, and $v\in p\inv(x)$.
It follows from \thmref{covering quotient}, \corref{skew universal}, and
Propositions~\ref{skew normalizer} and
\ref{orbit skew}
that $(\OO,p)$ is isomorphic to a relative skew product covering
$\LL\times_\e \pi(\LL,x)/p_*\pi(\OO,v)$.
\end{proof}

\subsection*{Gross-Tucker Theorem}

If the subgroup $H$ of $G$ is normal, then a relative skew product
$\LL\times_\e G/H$ may be regarded as an ordinary skew product associated to
the cocycle $\LL\to G/H$ obtained from $\e$ by composing with the quotient
homomorphism $G\to G/H$. In particular, with the notation from the above
proof, if the subgroup $p_*\pi(\OO,v)$ of $\pi(\LL,x)$ is normal, then the
given connected covering $(\OO,p)$ is isomorphic to a skew-product covering.
From \corref{covering normal}, we know that this will happen if and only if
$\aut(\OO,p)$ acts transitively on $p\inv(x)$. One situation where this is
obviously true is for a covering $\OO\to\OO/G$, where $G$ is a group acting
freely on a connected $k$-graph $\OO$.

While all this is a nice application of the general theory of connected
coverings, it cheats us out of the full truth: connectedness of $\OO$ is
unnecessary, as we'll show in the following result, a version of the
Gross-Tucker Theorem (for the graph version, see
\cite{gt}*{Theorem 2.2.2}).
In the disconnected case it is more efficient to give a
``bare-hands'' proof.
Actually, this result appears in \cite{kp:kgraph}*{Remark 5.6}, but
we prove it here
for the convenience of the reader,
since we will need this more general result for
$C^*$-coactions:

\begin{thm}
[Gross-Tucker Theorem]
\label{gross tucker}
Let $G$ be a group acting freely on a $k$-graph $\S$. Then the covering
$\S\to\S/G$ given by the quotient map
is isomorphic to a skew product covering $(\S/G)\times_\e
G\to\S/G$.
\end{thm}

\begin{proof}
The corresponding groupoid-action result, \lemref{cocycle action} below, is
easier, so here we merely indicate how the Gross-Tucker Theorem will follow.
Put $\LL=\S/G$, and
let $t\:\S\to \LL$ be the quotient map.
Then $G$ is a subgroup of $\aut(\S,t)$
acting freely and transitively on each set
$t\inv(x)$ for $x\in\LL^0$.
By the elementary \lemref{cocycle action} below,
the associated groupoid action
$(\S^0,\G(\LL))$ is isomorphic to a cocycle action
$(\LL^0\times_\e G,\G(\LL))$, so
the covering $t\:\S\to\LL$ is isomorphic to
a skew-product covering $\LL\times_\e G\to\LL$.
\end{proof}

We must pay the debt we incurred in the above proof:

\begin{lem}
\label{cocycle action}
Let $G$ be a group acting on the right of a groupoid action $(V,\G)$,
freely and
transitively on each set $V_x$ for $x\in\G^0$. Then $(V,\G)$ is isomorphic to
a cocycle action.
\end{lem}

\begin{proof}
We begin by choosing a cross-section of the map $V\to\G^0$:
for each $x\in\G^0$ pick $v_x\in V_x$.
Let $x,y\in\G^0$ and $a\in x\G y$.
Then both $av_y$ and $v_x$ are in $V_x$,
so by hypothesis there exists a unique
element $\e(a)$ of $G$ such that
$v_x\e(a)=av_y$.
We verify that the resulting map
$\e\:\G\to G$
is a cocycle:
if $x,y,z\in\G^0$, $a\in x\G y$, and $b\in y\G z$,
then
\begin{align*}
v_x\e(a)\e(b)
&=av_y\e(b)
=abv_z
=v_x\e(ab),
\end{align*}
so $\e(a)\e(b)=\e(ab)$ since $G$ acts freely.

Define $\p\:\G^0\times_\e G\to V$ by
$\p(x,g)=v_xg$.
To see that $\p$ is injective, let $(x,g),(y,h)\in\G^0\times G$, and assume
$\p(x,g)=\p(y,h)$. Then
\[
v_yhg\inv=v_x\in V_x.
\]
Since $v_y\in V_y$, so is $v_yhg\inv$. Thus we must have $x=y$, hence $g=h$
since $G$ acts freely.
To see that $\p$ is surjective,
let $x\in\G^0$ and $v\in V_x$. Since $G$ acts
transitively on $V_x$, we can choose $g\in G$ such that $v=v_xg$, and then
$v=\p(x,g)$.
For $x,y\in \G^0$,
$a\in x\G y$, and $g\in G$
we have
\begin{align*}
\p\bigl(a(y,g)\bigr)
=\p(x,\e(a)g)
=v_x\e(a)g
=av_yg
=a\p(y,g),
\end{align*}
so $\p$ intertwines the cocycle action and the given action.
Therefore $\p\:(\G^0\times_\e G,\G)\to(V,\G)$ is an isomorphism.
\end{proof}

\section{Coactions}
\label{sec:coactions}

In this section we discuss the implications of our results for the $C^
*$-algebras of $k$-graphs. There are by now several different classes of
$k$-graphs $\Lambda$ whose $C^*$-algebras $C^*(\Lambda)$ admit a satisfactory
structure theory, and our results apply to all of them. Indeed, of the usual
theory we need to know only that the core is AF and that the Gauge-Invariant
Uniqueness Theorem holds. Thus the results of this section apply to, in
increasing order of generality, the row-finite $k$-graphs without sources of
\cite{kp:kgraph}, the locally convex row-finite $k$-graphs of \cite{rsy}, and
the finitely aligned $k$-graphs of \cite{rsy_finitely}.

The main point of \cite{kqr:graph} is that a labelling of a graph
gives rise to a coaction on the graph $C^*$-algebra, and that moreover the
coaction crossed product is isomorphic to the $C^*$-algebra of the
skew-product graph. Here we adapt this to $k$-graphs.

For $C^*$-coactions we adopt the conventions of
\cites{qui:discrete,
ekq,
eq:induced,
qui:fullred}.
A \emph{coaction} of a group $G$ on a $C^*$-algebra $A$
is an injective nondegenerate homomorphism $\d$ of $A$ into the spatial
tensor product $A\otimes C^*(G)$
satisfying the \emph{coaction identity} $(\id\otimes \d_G)\d=(\d\otimes
\id)\d$, where $\d_G$ is the comultiplication on $C^*(G)$. For $g\in G$
the associated \emph{spectral subspace} of $A$ is
$A_g:=\{a\in A\mid \d(a)=a\otimes g\}$,
and the \emph{fixed point algebra} is $A^\d:=A_e$.
The disjoint union $\A:=\bigsqcup_{g\in G}A_g$ is a \emph{Fell bundle} in the
sense that $A_gA_h\subset A_{gh}$ and $A_g^*=A_{g\inv}$,
and the linear span $\Gamma_c(\A)$ is a dense $*$-subalgebra of $A$.
Define $\rho\:\A\to G$ by $\rho(a)=g$ if $a\in \A_g$.
The coaction $\d$ is called \emph{maximal} if the norm of $A$ is
the largest $C^*$-norm on the $*$-algebra $\Gamma_c(\A)$ \cite{ekq},
and \emph{normal} if $(\id\otimes \ll)\d$ is injective, where $\ll$
is the left regular representation of $G$ \cite{qui:fullred}.
For a subgroup $H$ of $G$, the Cartesian product $\A\times G/H$ is a Fell
bundle over the transformation
groupoid $G\times G/H$, with operations
\[
(a,\rho(b)gH)(b,gH)=(ab,gH)
\and
(a,gH)^*=(a^*,\rho(a)gH),
\]
and the linear span $\Gamma_c(\A\times G/H)$ is a $*$-algebra,
whose completion $A\times_{\d|}G/H$ in the largest $C^*$-norm is the
\emph{restricted crossed product} of $A$ by $\d$.
When $H=\{e\}$ the \emph{dual action} of $G$ on
the crossed product
$A\times_\d G$ is given by
$\what\d_h(a,g)=(a,gh\inv)$.

Let $\eta\:\LL\to G$ be a $k$-graph cocycle. The right
action of $G$ on $\LL\times_\e G$ discussed in
\propref{skew action}
induces an action $\gamma$ of 
$G$ on $C^*(\LL\times_\e G)$ such that
\[
\g_h\bigl(s_{(\ll,g)}\bigr)=s_{(\ll,gh\inv)}.
\]

\begin{thm}
\label{skew coaction}
Let $\e\:\LL\to G$ be a $k$-graph cocycle and $H$ a subgroup of $G$.
Then:
\begin{enumerate}
\item
there exists a unique coaction $\d=\d_\e$
of $G$ on $C^*(\LL)$ such that
\[
\d(s_\ll)=s_\ll\otimes\e(\ll)
\for
\ll\in\LL;
\]

\item
$C^*(\LL\times_\e G/H)\cong C^*(\LL)\times_{\d|} G/H$;

\item
if $H=\{e\}$ the above isomorphism is equivariant for the action
$\g$ of $G$ on $C^*(\LL\times_\e G)$ and the dual action $\what\d$ on
$C^*(\LL)\times_\d G$;

\item
the coaction $\d$ is both maximal and normal.
\end{enumerate}
\end{thm}

Our desired extension of \cite{dpr}*{Theorem 3.2} follows immediately from
\thmref{skew coaction} and \corref{connected skew}:

\begin{cor}
\label{covering crossed product}
Let $p\:\OO\to\LL$ be a connected covering, $x\in\LL^0$, and $v\in p\inv(x)$.
Then there exists a coaction $\d$ of $\pi(\LL,x)$ on $C^*(\LL)$ such that
\[
C^*(\OO)\cong C^*(\LL)\times_{\d|} \pi(\LL,x)/p_*\pi(\OO,v).
\]
\end{cor}

\begin{proof}[Proof of \thmref{skew coaction}]
(i)
It is routine to verify that
the assignment $\ll\mapsto s_\ll\otimes\e(\ll)$
gives a Cuntz-Krieger $\LL$-family
in $C^*(\LL)\otimes C^*(G)$, and hence determines a unique homomorphism
$\d\:C^*(\LL)\to C^*(\LL)\otimes C^*(G)$;
$\d$ is nondengenerate and satisfies
the coaction identity, and the
Gauge-Invariant Uniqueness Theorem shows that $\d$ is
injective.

(ii)
Define $\t\:\LL\times_\e G/H\to C^*(\LL)\times_{\d|} G/H$ by
\[
\t(\ll,gH)=(s_\ll,gH).
\]
It is routine to verify that this gives a Cuntz-Krieger
$(\LL\times_\e G/H)$-family in the $C^*$-algebra $C^*(\LL)\times_{\d|} G/H$,
hence determines a homomorphism
\[
\t\:C^*(\LL\times_\e G/H)\to C^*(\LL)\times_{\d|} G/H.
\]
The Gauge-Invariant Uniqueness Theorem shows that $\t$ is injective, and it is
obviously surjective.

(iii)
For the equivariance,
\[
\t\g_h\bigl(s_{(\ll,g)}\bigr)
=\t\bigl(s_{(\ll,gh\inv)}\bigr)
=(s_\ll,gh\inv)
=\what\d_h(s_\ll,g)
=\what\d_h\t\bigl(s_{(\ll,g)}\bigr).
\]

(iv)
We first show that the coaction $\d$ is maximal.
Let $\A$ be the Fell bundle associated to the coaction $\d$.
Since the spectral subspaces are linearly independent in $C^*(\LL)$,
$\Gamma_c(\A)$ sits inside $C^*(\LL)$ as a $*$-subalgebra, giving an obvious
representation of $\A$ in $C^*(\LL)$,
which in turn extends uniquely to a 
homomorphism $\pi\:C^*(\A)\to C^*(\LL)$.
For maximality it suffices,
by \cite{ekq}*{Proposition 4.2},
to show that
$\pi$ is injective.
The inclusion $\LL\hookrightarrow \A$ gives a map
$\rho_0\:\LL\to C^*(\A)$;
the image is a Cuntz-Krieger $\LL$-family,
because the Cuntz-Krieger relations
can be expressed within the
Fell bundle $\A$.
Thus there is a unique
homomorphism $\rho\:C^*(\LL)\to C^*(\A)$
such that $\rho(s_\ll)=\rho_0(\ll)$ for $\ll\in\LL$.
We show that $\pi$ is injective by showing that $\rho$ is a left inverse. By \lemref{spectral} below, the spectral subspaces of the coaction $\d$, hence the Fell bundle $C^*$-algebra, are generated by the image of $\rho_0$, so it suffices to observe that for all $\ll\in\LL$ we have
\[
\rho\circ\pi(\rho_0(\ll))=\rho(s_\ll)=\rho_0(\ll).
\]

The homomorphism $\pi:=(\id\otimes\ll)\d$ intertwines the gauge action $\a$
and the tensor-product action $\a\otimes\id$, and $\pi(p_v)\ne 0$ for every
vertex $v$, so the Gauge-Invariant Uniqueness Theorem implies that $\pi$ is
faithful. Thus $\d$ is normal.
\end{proof}

\begin{cor}
For the Fell bundle $\A$ of the coaction $\delta\:C^*(\Lambda)\to
C^*(\Lambda)\otimes C^*(G)$, we have $C^*(\A)=C^*_r(\A)$ \(so that $\A$ is
amenable in the sense of Exel \cite{exe:amenability}\).
\end{cor}

\begin{proof}
The maximality of $\delta$ says that $C^*(\A)=C^*(\Lambda)$, and the
normality that the regular representation $(\id\otimes\lambda)\circ \delta$ is
an isomorphism of $C^*(\Lambda)$ onto
$C^*_r(\A):=\ran((\id\otimes\lambda)\circ \delta)$.
\end{proof}

\subsection*{Decomposition}

We apply \thmref{skew coaction} to give an analogue for $k$-graphs
of Green's decomposition theorem \cite{gre:local}*{Proposition~1} in which
the subgroup need not be normal and no twist is required.

\begin{cor}
\label{decomposition}
Let $\e\:\LL\to G$ be a $k$-graph cocycle, $\d=\d_\e$ the associated
coaction of $G$ on $C^*(\LL)$, and
$H$ a subgroup of $G$.
Then there is a coaction $\ee$ of $H$ on
the restricted crossed product $C^*(\LL)\times_{\d|} G/H$ such that
\[
C^*(\LL)\times_{\d|} G/H\times_\ee H\cong C^*(\LL)\times_\d G,
\]
equivariantly for the dual action $\what\ee$ and the restricted dual
action $\what\d\big|H$.
\end{cor}

\begin{proof}
Since our aim is to apply \thmref{skew coaction}, we need a cocycle.
$H$ acts freely on $\LL\times_\e G$,
and we have
\[
(\LL\times_\e G)/H\cong \LL\times_\e G/H.
\]
Thus by the Gross-Tucker Theorem~\ref{gross tucker} (twice!)
there is a cocycle
$\k\:\LL\times_\e G/H\to H$ such that
\[
\LL\times_\e G/H\times_\k H\cong \LL\times_\e G.
\]
By \thmref{skew coaction},
letting $\d_\k$ denote the corresponding coaction of $H$ on $C^*(\LL\times_\e
G/H)$,
we have
\[
C^*(\LL\times_\e G/H)\times_{\d_\k} H\cong C^*(\LL)\times_\d G,
\]
equivariantly for $\what{\d_\k}$ and $\what\d\bigm|H$.
Appealing to the Gross-Tucker Theorem once more we have $C^*(\LL)\times_\e
G/H\cong C^*(\LL)\times_{\d|} G/H$; this isomorphism is equivariant for a
unique coaction $\ee$ of $H$ on $C^*(\LL)\times_{\d|} G/H$,
and the result follows.
\end{proof}

\begin{rem}
For $1$-graphs, Corollary~\ref{decomposition} reduces to
\cite{dpr}*{Corollary 3.6}, except that it uses the full rather than
reduced crossed product. Since \cite{dpr}*{Corollary 3.6} motivated a
general result for deompositions of crossed products by normal coactions
using the reduced crossed product \cite{dpr}*{Theorem 4.2}, it is tempting
to conjecture on the basis of Corollary~\ref{decomposition} that there is
a similar decomposition for maximal coactions using the full crossed product.
\end{rem}

The next corollary extends \cite{kp:kgraph}*{Theorem 5.7}:

\begin{cor}
With the above hypotheses, and $\g$ the action of $G$ on $C^*(\LL\times_\e
G)$ described before \thmref{skew coaction}, we have
\[
C^*(\LL\times_\e G)\times_{\g|}H
\cong C^*(\LL\times_\e G/H)\otimes \K(l^2(H)).
\]
\end{cor}

\begin{proof}
We have:
\begin{align*}
C^*(\LL\times_\e G)\times_{\g|}H
&\cong C^*(\LL)\times_\d G\times_{\what\d}H
\\&\cong (C^*(\LL)\times_{\d|} G/H)\times_\ee H\times_{\what\ee}H
\\&\cong \bigl(C^*(\LL)\times_{\d|} G/H\bigr)\otimes \K(l^2(H))
\\&\cong C^*(\LL\times_\e G/H)\otimes \K(l^2(H)),
\end{align*}
where we successively applied:
\thmref{skew coaction},
\corref{decomposition},
crossed-product duality,
and \thmref{skew coaction} again.
\end{proof}

\subsection*{Cohomology}

The theories of both graphs and groupoids (see, e.g.,
\cites{gt, pask-rao, ren:approach}), contain a notion of cohomology
of cocycles.
This is easily adapted to $k$-graphs,
and has ramifications for the associated coverings and coactions:
we call cocycles $\e,\kappa\:\LL\to G$ \emph{cohomologous} if there exists
a map
$x\mapsto\tau_x\:\LL^0\to G$ such that
\[
\tau_x\e(a)=\kappa(a)\tau_y
\forall a\in x\LL y.
\]
If we regard $\e$ and $\kappa$ as functors
then the map $\tau$ is
just a natural isomorphism from $\e$ to $\kappa$.
It is routine to verify that the map $(x,g)\mapsto (x,\tau_x g)$ gives a
$k$-graph
isomorphism $\LL\times_\e G\cong \LL\times_\k G$
which is equivariant for the associated actions of $G$, and
the unitary multiplier
$\sum_{x\in\LL^0}(x\otimes \tau_x)$
implements an exterior equivalence
between the associated coactions $\d_\e$ and $\d_\k$.

\subsection*{The gauge coaction}

We can view the degree functor
as a cocycle $d\:\LL\to\Z^k$. By \thmref{skew coaction} there is a unique
coaction
$\d=\d_d$ of
$\Z^k$ on $C^*(\LL)$ such that
\[
\d(s_\ll)=s_\ll\otimes d(\ll)
\for \ll\in\LL.
\]
We call $\d$ the \emph{gauge coaction} because the corresponding action of
$\mathbb{T}^k=\what{\Z^k}$ is the usual gauge action.

When $\LL$ is row-finite and has no sources, the following result is contained
in \cite{kp:kgraph}*{Theorem 5.5}.

\begin{thm}
\label{nuclear} Suppose that $\LL$ is a countable finitely aligned $k$-graph.
Then
$C^*(\LL)$ is nuclear,
and $C^*(\LL\times_d \Z^k)$ is AF.
\end{thm}

\begin{proof}
[Proof of nuclearity]
The fixed-point
algebra $C^*(\LL)^\d$ is the core, which
is  AF (see the proof of \cite{rsy_finitely}*{Theorem~3.1}).
Thus $C^*(\LL)^\d$ is in  particular nuclear, and
\cite{qui:discrete}*{Corollary 2.17} implies that $C^*(\LL)$ is also nuclear.
\pause

We will prove that $C^*(\LL\times_d \Z^k)$ is AF by proving that the isomorphic
algebra $C^*(\LL)\times_\d \Z^k$ is AF.
The proof
would not be hard if we had \emph{saturation} (see below
for the definition), for then the crossed product
would be Morita-Rieffel
equivalent to the fixed-point algebra.
However, in the general case we require a digression.

Recall from \cite{qui:discrete} that an
\emph{ideal property} is
a property $\P$ of $C^*$-algebras such that: (1) every $C^*$-algebra
has a largest ideal with $\P$, (2) $\P$ is inherited by ideals, and
(3) $\P$ is invariant under Morita-Rieffel equivalence.
The motivation for this definition was then, and remains for us here, that if
$\d$ is a coaction of a discrete group $G$ on a $C^*$-algebra $A$, then for
any ideal property $\P$, the crossed product $A\times_\d G$ has $\P$ if and
only if the fixed-point algebra $A^\d$ does.
It is shown in \cite{qui:discrete} that nuclearity is an ideal property,
and it is well-known that liminality and postliminality are ideal
properties.

\begin{prop}
\label{AF}
Among separable $C^*$-algebras, AF is an ideal property.
\end{prop}

\begin{proof}
For invariance under Morita-Rieffel equivalence, let $A\sim B$ with $A$
being AF. Since $A$ and $B$ are separable, we have
$A\otimes\K\cong B\otimes\K$. Since $A$ is AF, so is $A\otimes\K$, hence
$B\otimes\K$.
Thus the hereditary subalgebra $B$ is also AF,
by \cite{elliott}*{Theorem 3.1}.
The same result of Elliott shows that AF is
inherited by ideals.

We finish by showing that every $C^*$-algebra
$A$ has a largest AF
ideal, i.e., an AF ideal which contains every AF ideal. Claim:
if $I$ and $J$ are AF ideals of $A$, then the ideal $I+J$ is AF. Since
\[
(I+J)/I\cong I/(I\cap J),
\]
the quotient $(I+J)/I$ is AF. Thus
the extension
$I+J$ of
$(I+J)/I$
by
$I$
is AF, by results of Brown \cite{bro}
and Elliott \cite{elliott}.
We
pause to make this reference more
precise, because the required result must be pieced together:
Elliott proved in \cite{elliott}*{Corollary 3.3} that AF is closed
under extensions provided projections lift, and Brown proved that projections
do indeed lift (from an AF quotient by an AF ideal) --- actually, the full
proof of Brown's result is in \cite{effros}*{Section 9}.

Now let $I$ be the closed span of all AF ideals of $A$. Then $I$ is
certainly an ideal of $A$. By the above, $I$ is the closure of an
upward-directed union of AF ideals. Therefore $I$ is AF. By construction,
every AF ideal of $A$ is contained in $I$.
\end{proof}

We're now ready to prove that the crossed product by the gauge coaction is AF:

\resume{\thmref{nuclear}}
Since $C^*(\LL)^\d$ is AF
and $C^*(\LL)\times_\d\Z^k\cong C^*(\LL\times_d\Z^k)$,
the result follows from \cite{qui:discrete}*{Corollary 2.17},
because AF is an ideal property.
\end{proof}

When $k=1$,
the skew-product graph has no cycles,
so its $C^*$-algebra is AF by, for example,  \cite{dt}*{Corollary~2.13}.
If we had
a corresponding result concerning cycles for $k>1$,
this would give an alternate proof of
the second part of \thmref{nuclear} ---
but we don't.

\subsection*{Saturation}

Let $\d=\d_d$ be the gauge coaction of $\Z^k$ on $C^*(\LL)$.
Recall from \cite{qui:discrete} that
$\d$ is called \emph{saturated} if
\[
C^*(\LL)_{n+m}=\clspn C^*(\LL)_nC^*(\LL)_m
\righttext{for every}
n,m\in \Z^k,
\]
or equivalently if
\[
C^*(\LL)^\d=\clspn C^*(\LL)_nC^*(\LL)_n^*
\righttext{for every}
n\in \Z^k.
\]
If $\d$ is saturated then $C^*(\LL)^\d$
is Morita-Rieffel equivalent to $C^*(\LL)\times_\d\Z^k$ \cite{qui:discrete}.

Recall that a vertex $v$ of $\LL$ is called a \emph{source} if
$v\LL^n=\emptyset$ for some $n\in\N^k$,
and similarly a \emph{sink} if
some
$\LL^n v$ is empty.
The following result generalizes
\cite{kp}*{Proposition 2.8}.

\begin{prop}
\label{saturated}
Let $\d$ be the gauge coaction of $\Z^k$ on $C^*(\LL)$.
\begin{enumerate}
\item
If $\d$ is surjective, in particular
if $\LL$ has
either
no sources or no sinks,
every spectral subspace $C^*(\LL)_n$ for $n\in\Z^k$ is nontrivial.

\item
If
$\LL$ is row-finite and
has neither sources nor sinks, then $\d$ is saturated.
\end{enumerate}
\end{prop}

\begin{proof}
(i)
Let $n\in\Z^k$.
Choose $l\in\N^k$ with $n+l\ge 0$, then
$j\in\N^k$ with $j\ge n+l$ and $j\ge l$,
and then $\ll\in\LL^j$.
We can factor
\[
\ll=\m\n=\a\b
\midtext{with}
d(\n)=n+l
\and
d(\b)=l.
\]
Thus
\[
0\ne s_\ll s_\ll^*=s_\m s_\n s_\b^*s_\a^*,
\]
so that $s_\n s_\b^*$ is a nonzero element of
\[
C^*(\LL)_{n+l}C^*(\LL)_l^*
\subset C^*(\LL)_n.
\]

(ii)
Now assume that $\LL$
is row-finite and has neither sources nor sinks.
Let $l\in\Z^k$.
To see that $\d$ is saturated,
we must show that
\[
C^*(\LL)^\d\subset\clspn C^*(\LL)_lC^*(\LL)_l^*.
\]
By \lemref{spectral} below,
we have
\[
C^*(\LL)^\d=C^*(\LL)_0
=\clspn\{s_\ll s_\m^*\mid d(\ll)=d(\m), s(\ll)=s(\m)\},
\]
so it suffices to show that if
\[
d(\ll)=d(\m)=n
\midtext{and}
s(\ll)=s(\m)=v,
\]
then
$s_\ll s_\m^*\in \spn C^*(\LL)_l C^*(\LL)_l^*$.
Choose $m\in\N^k$ with $m\ge l$ and $m\ge n$.
Since $\LL$ is row-finite and has no sources,
\[
p_v=\sum_{\a\in v\LL^{m-n}}s_\a s_\a^*.
\]
Since $\LL$ has no sinks,
for each $\a\in v\LL^{m-n}$ we can choose $\n_\a\in \LL^{m-l}s(\a)$.
Then $p_{s(\a)}=s_{\n_\a}^*s_{\n_\a}$, so
\begin{align*}
s_\ll s_\m^*
&=s_\ll p_v s_\m^*
=\sum_{\a\in v\LL^{m-n}}
s_\ll s_\a s_{\n_\a}^* s_{\n_\a} s_\a^* s_\m^*
\\&\in \spn C^*(\LL)_n C^*(\LL)_{m-n} C^*(\LL)_{m-l}^*
C^*(\LL)_{m-l} C^*(\LL)_{m-n}^* C^*(\LL)_n^*
\\&\subset \spn C^*(\LL)_l C^*(\LL)_l^*.
\qedhere
\end{align*}
\end{proof}

In the above proof,
we applied the following characterization of spectral subspaces
for the gauge coaction:

\begin{lem}
\label{spectral}
Let $\e\:\LL\to G$ be a cocycle on the $k$-graph $\LL$, and let $\d=\d_\e$ be
the associated coaction of $G$ on $C^*(\LL)$.
Then for all $g\in G$,
\[
C^*(\LL)_g
=\clspn\{s_\ll s^*_\m\mid \e(\ll)\e(\m)\inv=g\}.
\]
\end{lem}

\begin{proof}
Obviously a product $s_\ll s^*_\m$ is in $C^*(\LL)_g$ if and only if
$\e(\ll)\e(\m)\inv=g$, so the left hand side contains the right.

Recall from \cite{qui:discrete} that
there is a
bounded linear projection
\[
E_g=(\id\otimes\Chi_{\{g\}})\d\:C^*(\LL)\to C^*(\LL)_g,
\]
where here
the characteristic function
$\Chi_{\{g\}}$ is regarded as a linear functional on
$C^*(G)$.
Any $a\in C^*(\LL)_g$ can be approximated by a linear combination
$\sum_1^nc_is_{\ll_i}s_{\m_i}^*$, and then
\begin{align*}
a
&=E_g(a)
\approx E_g\biggl(\sum_1^nc_is_{\ll_i}s_{\m_i}^*\biggr)
=\sum_1^nc_iE_g(s_{\ll_i}s_{\m_i}^*)
\\&=\sum\{c_is_{\ll_i}s_{\m_i}^*\mid \e(\ll_i)\e(\m_i)\inv=g\},
\end{align*}
which is in the right hand side. The result follows.
\end{proof}

\begin{ex}
Row-finiteness is necessary in \propref{saturated}~(ii):
let $\LL$ be the 1-graph
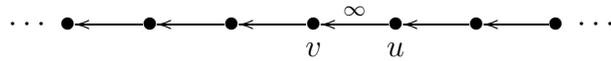
\[
\xymatrix{
& \point \ar@{}[l]|{\displaystyle{\cdots}}
& \point \ar[l]
& \point \ar[l]
& \point \ar[l] \ar@{}[d]|(.3){\displaystyle{v}}
& \point \ar[l]_\infty \ar@{}[d]|(.3){\displaystyle{u}}
& \point \ar[l]
& \point \ar[l]
& \ar@{}[l]|{\displaystyle{\cdots}}
\\
&&&&&
}
\]
in which there are infinitely many edges from $u$ to $v$ (and the graph
extends indefinitely to the right and left).
The projection $p_v$ is in the fixed-point algebra $C^*(\LL)^\d$, and cannot
be approximated in norm by a linear combination of products $s_es_f^*$ for
edges $e,f\in v\LL u$. It follows that
\[
C^*(\LL)^\d\ne \clspn C^*(\LL)_1C^*(\LL)_1^*,
\]
so $\d$ is not saturated.
\end{ex}



\end{document}